\documentclass[12pt,leqno]{amsart}
\usepackage{amssymb,amsmath}

\begin{document}

\title[communication networks and Hilbert
modular Forms]{Designing communication networks via
Hilbert modular Forms\thanks{The author is
partially supported by a Binational Israel-USA
Foundation grant}}

\author[R. Livn\'e]{Ron Livn\'e}
 \address{Mathematics Institute\\
The Hebrew University of Jerusalem\\
Givat-Ram, Jerusalem  91904\\  Israel}

\begin{abstract}
Ramanujan graphs, introduced by Lubotzky, Phillips
and Sarnak, allow the design of efficient
communication networks. In joint work with B.
Jordan we gave a higher-dimensional generalization.
Here we explain how one could use this
generalization to construct efficient communication
networks which allow for a number of verification
protocols and for the distribution of information
along several channels. The efficiency of our
network hinges on the Ramanujan-Petersson
conjecture for certain Hilbert modular forms. We
obtain this conjecture in sufficient generality to
apply to some particularly appealing constructions,
which were not accessible before.
\end{abstract}

\maketitle

\newcounter{soussec}
\newcommand{\soussection}[1]{\par \smallskip
\par \noindent {\bf \arabic{section}.\arabic{soussec}\
#1 \hspace*{1em}} \refstepcounter{soussec}}

\newtheorem{theorem}{Theorem}[section]
\newtheorem{proposition}[theorem]{Proposition}
\newtheorem{variant}[theorem]{A variant}
\newtheorem{fact}[equation]{}
\newtheorem{corollary}[theorem]{Corollary}
\newtheorem{definition}[theorem]{Definition}
\newtheorem{example}[theorem]{Example}
\newtheorem{examples}[theorem]{Examples}
\newtheorem{conditions}[theorem]{Conditions}
\newtheorem{remark}[theorem]{Remark}
\newtheorem{remarks}[theorem]{Remarks}
\newtheorem{lemma}[theorem]{Lemma}
\newtheorem{maintheorem}[theorem]{Main Theorem}

\newenvironment{pf}{\noindent {\it Proof:}\/ }
{\par \par \medskip\par }

\newcommand{\ra}{\rightarrow}
\newcommand{\sra}[1]{\stackrel{#1}{\rightarrow}}
\newcommand{\la}{\leftarrow}
\newcommand{\lra}{\longrightarrow}
\newcommand{\ov}{\overline}
\newcommand{\bs}{\backslash}
\newcommand{\matr}[4]
{\left[\begin{array}{cc} #1 & #2 \\ #3 & #4
              \end{array} \right]}
\newcommand{\smatr}[4]{\Bigl[\begin{array}
{@{\hspace{0.1em}}c@{\hspace{0.5em}}c@{\hspace{0.1em}}}
   \scriptstyle #1 & \scriptstyle #2 \\
   \scriptstyle #3 & \scriptstyle #4 \end{array}
   \Bigr]}
\newcommand{\set}[1]{\{\,#1\,\}}
\newcommand{\ip}[2]{{\langle #1 ,#2\rangle}}

\newcommand{\Aut}{\operatorname{Aut}}
\newcommand{\Disc}{\operatorname{Disc}}
\newcommand{\edg}{\operatorname{Ed}}
\newcommand{\edgo}[1]
        {\operatorname{Ed}^{\rm o}_{#1}}

\newcommand{\Gal}{\operatorname{Gal}}
\newcommand{\GL}{\operatorname{GL}}
\newcommand{\Gr}{\operatorname{Gr}}
\newcommand{\Hom}{\operatorname{Hom}}
\newcommand{\Ima}{\operatorname{Im}}
\newcommand{\Ker}{\operatorname{Ker}}
\newcommand{\Lef}{\operatorname{Lef}}
\newcommand{\Ltwo}{\operatorname{L}_2}
\newcommand{\Mat}{\operatorname{Mat}}
\newcommand{\Matwo}{\Mat_{2\times 2}}
\newcommand{\Nm}{\operatorname{Nm}}
\newcommand{\PGL}{\operatorname{PGL}}
\newcommand{\PSL}{\operatorname{PSL}}
\newcommand{\Res}{\operatorname{Res}}
\newcommand{\res}{\operatorname{res}}
\newcommand{\SL}{{\operatorname{SL}}}
\newcommand{\St}{\operatorname{St}}
\newcommand{\SU}{\operatorname{SU}}
\newcommand{\Symm}{\operatorname{Symm}}
\newcommand{\Tr}{\operatorname{Tr}}
\newcommand{\Ver}{\operatorname{Ver}\,}
\newcommand{\sign}{\operatorname{sign}}
\newcommand{\FGal}{F^{{\rm Gal}}}

\newcommand{\pard}[1]{\partial_{#1}}
\newcommand{\parad}[1]{\partial^\ast_{#1}}
\newcommand{\der}{{d}}
\newcommand{\derad}{\der^\ast}
\newcommand{\parj}{\pard{j}}
\newcommand{\paradj}{\parad{j}}
\newcommand{\invlim}
    {\underleftarrow{\vphantom{_U}\lim} \vphantom{\lim}}
\newcommand{\ran}{{\rm an}}
\newcommand{\rf}{{\rm f}}
\newcommand{\Const}{{\rm C}}
\newcommand{\vk}{\vec{k}}
\newcommand{\vs}{\vec{s}}
\newcommand{\Fpb}{\ov{\FF}\vphantom{\FF}_p}
\newcommand{\Qpb}{\ov{\QQ}\vphantom{\QQ}_p}

\newcommand{\Id}{{\rm Id}}
\newcommand{\Lap}[1]{\square_{#1}}
\newcommand{\Lapj}{\Lap{j}}
\newcommand{\Lapt}[1]{\Lap{{\rm tot},#1}}
\newcommand{\ocell}[1]{{\Sigma_{#1}}}
\newcommand{\ZtwoI}{(\ZZ/2\ZZ)^I}
\newcommand{\inv}[1]{{\rm inv}_{#1}}
\newcommand{\topf}[1]{{\rm top}_{#1}}
\newcommand{\tree}{{\bf T}}
\newcommand{\botf}[1]{{\rm bot}_{#1}}
\newcommand{\invj}{\inv{j}}
\newcommand{\topj}{\topf{j}}
\newcommand{\botj}{\botf{j}}
\newcommand{\ktype}{$(k_1,\ldots,k_g)$\/}
\newcommand{\HQ}{{\rm H}}
\newcommand{\pib}{{\ov{\pi}}}
\newcommand{\nub}{{\ov{\nu}}}
\newcommand{\gammab}{{\ov{\gamma}}}
\newcommand{\SN}{{\SL_2(\FF_N)}}
\newcommand{\PSN}{{\PSL_2(\FF_N)}}
\newcommand{\phib}{{\ov{\phi}}}
\newcommand{\cMF}{{\cM_{0,F}}}
\newcommand{\leg}[2]{{\left(\frac{#1}{#2}\right)}}
\newcommand{\symsp}{\cD}
\newcommand{\HB}{\cX}
\newcommand{\Frob}{{\rm Frob}}
\newcommand{\triv}[1]{{{\mathcal T}_{#1}}}
\newcommand{\Gm}{\GG_m}

\newcommand{\calA}{{\mathcal A}}
\newcommand{\calB}{{\mathcal B}}
\newcommand{\calD}{{\mathcal D}}
\newcommand{\calH}{{\mathcal H}}
\newcommand{\calI}{{\mathcal I}}
\newcommand{\calJ}{{\mathcal J}}
\newcommand{\calL}{{\mathcal L}}
\newcommand{\calM}{{\mathcal M}}
\newcommand{\calO}{{\mathcal O}}
\newcommand{\calT}{{\mathcal T}}
\newcommand{\calX}{{\mathcal X}}
\newcommand{\cA}{{\calA}}
\newcommand{\cB}{{\calB}}
\newcommand{\cD}{{\calD}}
\newcommand{\cH}{{\calH}}
\newcommand{\cI}{{\calI}}
\newcommand{\cL}{{\calL}}
\newcommand{\cM}{{\calM}}
\newcommand{\cO}{{\calO}}
\newcommand{\cT}{{\calT}}
\newcommand{\cX}{{\calX}}
\newcommand{\bG}{{\bf G}}
\newcommand{\bH}{{\bf H}}
\newcommand{\bV}{{\bf V}}
\newcommand{\bZ}{{\bf Z}}

\renewcommand{\AA}{{\mathbb A}}
\newcommand{\CC}{{\mathbb C}}
\newcommand{\FF}{{\mathbb F}}
\newcommand{\GG}{{\mathbb G}}
\newcommand{\HH}{{\mathbb H}}
\newcommand{\PP}{{\mathbb P}}
\newcommand{\QQ}{{\mathbb Q}}
\newcommand{\RR}{{\mathbb R}}
\newcommand{\TT}{{\mathbb T}}
\newcommand{\ZZ}{{\mathbb Z}}

\newcommand{\sig}{\sigma}
\newcommand{\lam}{\lambda}
\newcommand{\vpi}{\varpi}
\newcommand{\FFq}{\FF_q}

\newcommand{\qn}{{\hat 1}}
\newcommand{\qi}{{\hat\imath}}
\newcommand{\qj}{{\hat\jmath}}
\newcommand{\qk}{{\hat k}}


The concept of a Ramanujan graph was introduced and
studied by Lubotzky, Phillips and Sarnak (LPS) in
\cite{LPS}. These are $r$\/-regular graphs for
which the nontrivial eigenvalues $\lambda\neq \pm
k$ of the adjacency matrix satisfy the bounds
$|\lambda| \leq 2\sqrt{r-1}$. In many aspects these
bounds are optimal and natural. For example, the
adjacency matrix is the combinatorial analog of the
Laplace operator, and the bounds parallel the
(conjectured) Selberg bounds for the Laplacian on
Riemann surfaces. The main result of LPS was an
explicit construction of $p+1$\/-regular such
graphs, $p\equiv 1\mod{4}$ a prime, through the
arithmetic of quaternion algebras over the rational
numbers. From the point of view of Communication
Network Theory, the arithmetic examples are
particularly interesting: all Ramanujan graphs are
super-expanders; but in addition the examples have
many other useful properties, for example very good
expansion constants and large girth. Thus they can
be used to design efficient communication networks.

The Ramanujan property for the LPS examples hinges
on the truth of the Ramanujan-Petersson conjecture
for an appropriate space of modular forms of weight
$2$ over $\QQ$. Let $f$ be a weight $2$ holomorphic
cuspidal Hecke eigenform on a congruence subgroup.
The conjecture is that for any prime $p$ not
dividing the level the Hecke eigenvalue $a_p$
satisfies $|a_p| \leq 2 \sqrt{p}$\/. Eichler and
Shimura reduced the conjecture to Weil's results on
the absolute value of Frobenius eigenvalues for
curves over a finite fields. They achieved this for
a fixed form and all but a finite, unspecified set
of primes $p$. The proof requires a deep study of
the reduction modulo $p$ of modular curves and
Hecke correspondences. Igusa subsequently showed
that the method applied for all primes not dividing
the level of the form (see \cite{Del} for a more
modern exposition and a generalization to forms
over $\QQ$ of any weight). Igusa's part is
essential for the applications to communication
networks because there one first chooses the prime
$p$ and only then an increasing sequence of levels
prime to it. This gives an infinite family of
Ramanujan graphs of {\em fixed} regularity $p+1$.

In a joint work with B. Jordan (\cite{JL8}) we gave
a higher dimensional generalization of this theory
to $(r_1,\dots,r_g)$\/-regular cubical complexes.
These are cell complexes locally isomorphic to a
product of $r_i$\/-regular trees. Here there are
partial Laplacians, one for each tree factor. The
parallel bounds for their eigenvalues define the
notion of being Ramanujan. We then constructed
infinite towers of explicit arithmetic examples
whenever each $r_i$ is of the form $q_i + 1$ for a
prime power $q_i$. We used quaternion algebras over
totally real fields and hence our generalization
required the Ramanujan-Petersson conjecture for
Hilbert modular forms of multi-weight
$(2,\dots,2)$. In place of the work of Eichler,
Shimura and Igusa we invoked analogous results of
Carayol (\cite{Car}) on the reduction modulo $p$ of
Shimura curves over totally real number fields,
which as before allow to use Weil's results.

Carayol's method imposes on the form a technical
assumption which forced us to exclude certain
natural and particularly appealing examples.
Building on ideas of Langlands, Brylinski and
Labesse (\cite{BL}) obtained results which are free
from such assumptions, by reducing the conjecture
to Deligne's bounds for Frobenius eigenvalues on
the intersection homology of Hilbert-Blumenthal
varieties.
However in their proof they used in an essential
way the Satake-Baily-Borel compactification which
was available only over $\QQ$. This forced them to
exclude for any fixed form an unspecified finite
set of primes, which makes their results
inapplicable to us as was explained above. One of
our aims here is to point out that a modification
of these arguments proves the following
\begin{theorem}
\label{RPp}
 Let $F$ be a totally real field of degree $d$ over $\QQ$. Let
$f$ be a holomorphic Hilbert cuspidal Hecke
eigenform over $F$ of multi-weight
$(k_1,\dots,k_d)$, where the $k_i$\/'s are integers
$\geq 2$ and all of the same parity. Let $v$ be a
prime of $F$ of degree $d_v=[F_v:\QQ_p]$ whose
residual characteristic $p$ is prime to the
discriminant $\Disc F$ of $F$ and to the level $N =
N(f)$ of $f$. Let $\lam_v(f)$ denote the eigenvalue
of the Hecke operator $T_v$ belonging to $f$. Then
the Ramanujan-Petersson conjecture holds for $f$
and $v$, namely we have $|\lam_v(f)| \leq
2p^{(k-1)d_v/2}$, where $k=\max_i k_i$.
\end{theorem}

For the theorem to be meaningful we must specify
the normalizations made in defining $T_v$ and
$\lam_v$. However Hilbert modular forms are best
approached through representation theory, and we
have in fact opted to define $T_v$ and $\lam_v$
only representation-theoretically. This spares us
the tedious task of making ``classical''
definitions and then comparing them to the
representation-theoretic ones. The version of the
theorem we will actually prove is therefore
Theorem~\ref{RP} below.

Here is an outline of the article. In Section 1 we
discuss the significance of the higher dimensional
theory of \cite{JL8} to communication networks,
especially in the two-dimensional case. In Section
2 we prove the above theorem. This enables us to
handle in Section~\ref{secspec} a particularly
pretty example which was left unsettled in
\cite{JL8}. The result is a mixture of graph
theory, automorphic forms, algebraic geometry, and
the arithmetic of quaternions over number fields.
In our attempt to make this work accessible to a
mixed audience and yet not too lengthy we have
undoubtedly failed to give the right level of
detail to any particular reader. We can only ask
for indulgence in this matter.

The debt this work owes to my long term
collaboration with B. Jordan, in particular to
\cite{JL8}, should be evident. I started thinking
about this problem at the instigation of P. Sarnak.
Different approaches to the problem were
subsequently suggested by D. Blasius and C.-L. Chai
(independently). P. Deligne suggested to use
compactly supported cohomology. Conversations with
them and with J. Bernstein, G. Faltings, G. Harder,
and N. Katz were encouraging and helpful. It is a
pleasure to thank them all.

\section{Cubical complexes and communication networks}
A basic problem in communication networks is to
design explicit {\em super-expanders}. For a given,
arbitrarily large set of nodes one wants to connect
each node to a fixed number $r$ of ``neighbor''
nodes, so that information can spread fast over the
resulting network $N$. A commonly used quantitative
measure of efficiency is the expansion constant ---
the largest real constant $c>0$ so that each subset
$A$ of the nodes of $N$ having $|A| \leq |N|/2$
nodes has at least $c |A|$ (new) neighbors. One
seeks sequences of networks $N_k$ where
$|N_k|\ra\infty$ with $k$, while $c$ is independent
of $k$ and is as large as possible. Even though one
can define the meaning of a {\em random} network on
$n$ vertices and prove that random networks are
good expanders, it is quite hard to get explicit
examples. The problem whether a random network is
Ramanujan is open, although it is known to be
``almost Ramanujan'' as $r\ra \infty$. More
precisely, $\lambda \leq 2\sqrt{r-1} + \log(r-1) +
{\rm Const}$ with probability $\ra 1$ with $N$ (see
\cite{Fri}).

In this work we shall explain how to design any
number $g \geq 1$ of efficient communication
networks on the same arbitrarily large set $N_k$ of
nodes. Thinking of the connections of each network
as having a different color, each node will have a
fixed number $r_i$ of color $i$ neighbors. Each
color will be a Ramanujan network, and hence a
super-expander. But in addition, any set of $h$
edges of different colors starting at the same
vertex can be completed to an $h$\/-cube. For
example, given a red edge $e_r$ and a blue edge
$e_b$ emanating from the same vertex, there will be
a unique pair of red and blue edges $e'_r$, $e'_b$
completing them to a square. In other words, the
origins $o(.)$ and ends $t(.)$ of the edges will
satisfy $o(e'_b) = t(e_r)$, $o(e'_r) = (e_b)$, and
$t(e'_b) = t(e'_r)$. We say that the network
satisfy the cube property (or the square property
for $2$ colors). Like the LPS examples, the ones we
give have good expansion properties and large
girth.

In practical applications the cube property permits
a variety of potential applications
\begin{enumerate}
\item To send information along (say) the red network, and
verify, using a hash function, a compatibility
condition via the blue channel.
\item To send divide the information to $g$ parts and
send each part individually. One could arrange so
that only the combination of the various elements
would lead to a meaningful whole.
\end{enumerate}

It seems a difficult problem to construct by random
methods systems of $(k_1,\dots,k_g)$\/-regular
networks satisfying the cube condition or to
parametrize the space of such systems.  The key
point is to realize that the cube condition makes
natural the introduction of higher dimensional
cubical complexes into the situation. Let $r_i \geq
3$ be a sequence of length $g$ of integers, let
$\tree_{r_i}$ be a regular tree of regularity
$r_i$, and set $\tree = \prod_i\tree_{r_i}$. Then
an $(r_1,\dots,r_g)$\/-regular cubical complex in
the sense of \cite{JL8} is a complex $X$ in which
each connected component is isomorphic to a
quotient $\Gamma\bs\tree$ for a discrete, torsion
free subgroup of $\prod_i \Aut\tree_{r_i}$. For
simplicity we also require the {\em parity
condition}: the $i$\/th component of any element of
$\Gamma$ moves each vertex of $\tree_{r_i}$ to an
even distance.  The $1$\/-skeleton of $X_\cI$ is a
collection of $g$ graphs, or communication
networks, on the same set of vertices. The parity
assumption makes these graphs bipartite. The
vertices have a ``multiparity'', an element of
$\set{0,1}^g$; the $i$\/th graph is
$r_i$\/-regular, and its $2^{g-1}$ connected
components consist of the vertices in which all
parities except for the $i$\/th have been fixed.
Its edges come from the $i$\/th tree factor of
$\tree$. Most importantly, the $g$ graphs satisfy
the cube property.

Theorem 3.1 of loc.cit.\ provides us with the
following general class of arithmetic examples:
\begin{theorem}
Let $B$ be a totally definite quaternion algebra
over a totally real field $F$, and let $S =
\set{v_1,\dots,v_g}$ be a nonempty set of $g$
distinct finite primes of of $F$ so that $B$ is
split at each $v_i$. Let $q_i$ be the cardinality
of the residue field of $F_{v_i}$ and set $r_i =
q_i + 1$. let $\Gamma$ be an $S$\/- arithmetic
torsion-free congruence subgroup of the algebraic
group over $\QQ$ of the norm $1$ elements in $B$.
Then $\Gamma$ acts on $\tree = \prod_i \tree_{r_i}$
and the resulting quotient $\Gamma\bs\tree$ is an
irreducible $(r_1,\dots,r_g)$\/-regular complex
with parities. Moreover if an appropriate space of
Hilbert cusp forms (of weight $(2,\dots,2)$)
satisfies the Ramanujan-Petersson conjecture, then
$\Gamma\bs\tree$ is Ramanujan.
\end{theorem}

The assumptions on $B$ mean that $B\otimes F_v$ is
isomorphic to the Hamilton quaternion algebra $\HH$
for any infinite prime $v$ of $F$, and that $B
\otimes F_{v_i}$ is isomorphic to $\GL(2,F_{v_i})$.
Let $\cO$ be the order of elements of $F$ integral
away from $v_1$, \dots, $v_g$ and let $\cM$ be a
maximal $\cO$\/-order of $B$. Then $\Gamma$ is
taken to be a sufficiently small congruence
subgroup of the group of norm $1$ elements of
$\cM$.

The description of our complexes as quotients of
infinite complexes by infinite groups can be
replaced by a finite description, which is more
elementary and very convenient for explicit
calculation. This finite description takes a
particularly simple shape if certain assumptions
are made (see \cite{JL8} for the details): let us
assume there is an ideal $N_0\neq 0$ of $\cO_F$,
prime to the $v_j$\/'s (we allow $N_0 = \cO_F$),
such that the following holds:
\begin{conditions}
\label{specase}
\begin{enumerate}
\item Every ideal of $F$ has a totally positive
generator $\equiv 1\mod N_0$.
\item The class number of $B$ is $1$.
\item The units $\cM^\times$ of a maximal
order $\cM$ of $B$ map onto $(\cM/N_0\cM)^\times$,
with the kernel being contained in the center
$\cO_F^\times$ of $\cM^\times$.
\end{enumerate}
\end{conditions}

Fix a totally positive generator $\pi_j$ which is
$\equiv 1\mod{N_0}$ for each prime ideal $v_j$. We
then have the following:
\begin{proposition}
\label{generators} {\rm 1.}\ For each $1\leq j \leq
g$ there are exactly $r_j$ {\rm (}principal\/{\rm
)} ideals $P_{j,i}$, $1\leq i \leq r_j$ of $\cM$
with norm $v_j$. We can moreover choose generators
\mbox{$\vpi_{j,i}\equiv 1 \mod{N_0\cM}$} for
$P_{j,i}$ whose norm is $\pi_j$.

\noindent {\rm 2.}\ For every permutation $\sigma$
of $\set{1,\dots,g}$ and any sequence of indices
$i_1,\dots,i_g$, with $1\leq i_j \leq r_j$, there
is a {\rm (}unique\/{\rm )} sequence $i'_1,\dots,
i'_g$, with $1\leq i'_j \leq r_j$, and a {\rm
(}unique\/{\rm )} unit $u\in \cO_F^\times$,
satisfying $u\equiv 1\mod{N_0}$, so that we have
the following arithmetic counterpart of the cube
condition:
\begin{equation}
\label{permut}
 \vpi_{\sig(1),i_1}\cdots \vpi_{\sig(g),i_g} =
         u\vpi_{1,i'_1}\cdots \vpi_{g,i'_g}\,.
\end{equation}
\end{proposition}

Let $N_1$ be a prime ideal of $\cO_F$ prime to the
$v_j$\/'s and to $N_0$, and set $N = N_0N_1$. Let
$A$ be the subgroup of $(\cO_F/N_1)^\times$
generated by the images modulo $N_1$ of the
$\pi_j$\/'s. Let $B$ be the subgroup of scalars in
$(\cM/N_1\cM)^\times$ generated by the images
modulo $N_1\cM$ of the $\pi_j$\/'s and by those
units of $\cO_F$ which are congruent to $1$ modulo
$N_0$. Set $H = \set{g\in (\cM/N_1\cM)^\times \mid
\Nm(g) \in A}
       / B$.
It is isomorphic to one of the groups
$\SL_2(\cO_F/N_1)$, $\PSL_2(\cO_F/N_1)$, or
$\PGL_2(\cO_F/N_1)$. (The examples of
Section~\ref{secspec} are of the $\PSL_2$ type.) We
now have the following

\begin{proposition}
\label{explicit} {\rm 1}.\ The vertices of $X(N)$
are the elements of $H$.

\noindent {\rm 2}.\ The {\rm (}oriented\/{\rm )}
edges of direction $j$ of $X(N)$ are the pairs
$(v,i_j)$, where $v\in H = \Ver X(N)$ and  $1\leq
i_j \leq r_j$.
\end{proposition}

In other words, each of the graphs is the Cayley
graph for {\em the same} group $H$ for a different
set of generators. This description is particularly
convenient for explicit calculation. Let us
emphasize that finite description exists even when
condition (\ref{specase}) does not hold (see
\cite{JL6} for the form it takes in the
prototypical one-dimensional case), but it is
messier. By \cite[Theorem~3.1(4)]{JL8} the
resulting graphs $\Gr_i(X)$ are Ramanujan provided
that appropriate spaces of cusp forms satisfy the
Ramanujan-Petersson conjecture. The statement that
fixing the parities gives connectedness follows
from Theorem 3.1(1) there.

\section{Hilbert modular forms}
\label{secHMF}
 \soussection{Hilbert-Blumenthal schemes} For
the generalities on Hilbert-Blumenthal schemes
which follow see \cite{DP,Rap}. Let $\cH_\pm =
\PP^1(\CC) - \PP^1(\RR)$ be the union of the upper
and the lower half complex planes. Let $\Sigma$ be
the set of the $d$ embeddings of $F$ into a Galois
closure $\FGal$, itself considered as a subfield of
$\RR$. Denote by $\bG$ the group scheme
$\Res_{\cO_F/\bZ} \GL_2$ over $\bZ$. Then
$\bG(\RR)$ is isomorphic to $\GL_2(\RR)^\Sigma$,
and $\bG(\RR)$ acts on $\symsp = (\cH_\pm)^\Sigma$
componentwise through the resulting $d$ M\"obius
transformations. Let $U$ be a compact open subgroup
of the finite ad\`eles $\bG(\AA_\rf) =
\GL_2(\AA_\rf \otimes F)$. Then the
Hilbert-Blumenthal complex space $\HB_U^\ran =
\bG(\QQ)\bs (\bG(\AA_\rf)/U \times \symsp)$ is
nonsingular if $U$ is sufficiently small. It has a
natural structure of a quasi-projective variety
which admits a canonical model over $\QQ$ in the
sense of Shimura (\cite{Del1}). In fact, if $U$
contains the principal congruence subgroup of level
$N$ and $D = \Disc F$ then there is a
moduli-theoretic interpretation of $\HB_U^\ran$ as
a coarse moduli space, parametrizing principally
polarized abelian $d$\/-folds with an
$\cO_F$\/-action and some level $N$ structure. This
yields a model $\HB_U$ over $\ZZ[1/ND]$ so that
$\HB_U \times \QQ$ is Shimura's canonical model. If
$U$ is sufficiently small then $\HB_U$ is a fine
moduli scheme, smooth and quasi-projective over
$\ZZ[1/ND]$. The Hecke algebra $\TT_U$ of
distributions on $U\bs \bG(\AA_\rf)/U$ acts on
$\HB_U\times \QQ$ through correspondences. We shall
suppose that $U$ is a product $U_pU^p$ of a part
$U_p\subset \bG(\QQ_pp$ at $p$ and a part
$U^p\subset \bG(\AA_f^p)$ away from $p$. Then an
element of $\TT_U$ coming from
$U_p\bs\bG(\QQ_p)/U_p$ (which we will call a Hecke
operator at $p$) acts on $\HB_U \times \ZZ[1/NDp]$
through its standard moduli interpretation.

Now let $p$ be a rational prime not dividing $ND$.
Then $\HB_U$ has good reduction modulo $p$. In
\cite{Lan} Langlands computed the number of fixed
points $T\times \Frob_p$ on $\HB_U \times \Fpb$ for
any Hecke operator $T$ away from $p$, i.e. $T$
coming from $U^p\bs\bG(\AA_f^p)/U^p$. In fact
Langlands allows more general groups than $\GL_2$
and he also allows algebraic local system on
$\HB_U$. This second generalization enables one to
handle cusp forms of weight $(k_1,\dots,k_d)$
rather than $(2,\dots,2)$. Even though the trivial
local system in the $\GL_2$ case is all that we
need here, there is no advantage in restricting to
it, and it has a potential use for the Ramanujan
local systems of \cite{JL8}. On the other hand we
decided to restrict to the $\GL_2$ case so as to
simplify our exposition by avoiding all mention of
$L$\/-packets and endoscopy. The interested reader
can see (\cite{BL, Lan}) that this restriction is
unnecessary. Thus let $\xi$ be an irreducible
algebraic representation of $\bG$ (all are defined
over the Galois closure $\FGal$ of $F$). We will
always assume that $U$ is sufficiently small and
that $\xi$ is trivial on the (Zariski closure) of
the units $\bZ(\cO_F)$ in the center $\bZ$ of
$\bG$. Then $\xi$ defines a local system
\[\bV_\xi^\ran = \bG(\QQ)\bs (\bG(\AA_\rf)/U \times \symsp
\times \xi)\]
 on $\HB_U^\ran$. The possible $\xi$\/'s can be described
explicitly: identify $\bG \times F \simeq \GL_2^d$
and let $\Symm^k$ be the representation of $\GL_2$
on the polynomials of degree $k\geq 0$ in two
variables. Then $\xi$ is a tensor product
$\otimes_{i=1}^d \xi_i$ of irreducible
representations $\xi$ of the $\GL_2$ factors, and
$\xi_i \simeq \Symm^{k_i-2} \otimes
\det^{\alpha_i}$, with $k_i \geq 2$. The condition
of being trivial on the units of $\bZ$ translates
into $k_i -2 + 2\alpha_i = \Const$ with $C$ a
constant. The $k_i$\/'s must all have the same
parity as $\Const$, and $\alpha_i = (\Const +2 -
k_i)/2$. Observe that the central character of
$\xi$ maps $t = (t_1,\dots,t_d)$ to
$\Nm(t)^\Const$, where $\Nm(t) = t_1\dots t_d$.

The $\xi$\/'s above have a moduli theoretic
interpretation and hence $\ell$\/-adic analogs
$\bV_{\xi,\lam}$, which are lisse sheaves of
$2$\/-dimensional $\FGal_\lambda$ vector spaces
over $\HB_U[1/\ell]$. Here $\lambda$ is any prime
of $\FGal$ lying above $\ell$. For this let $f:\cA
\ra \HB = \HB_U[1/\ell]$ be the universal abelian
variety (which we may assume exists by adding
auxiliary level structure). Then $R^1f_*\QQ_\ell$
is a lisse sheaf of $F\otimes\QQ_\ell$\/-modules of
rank $2$. For every $\sigma\in\Sigma$ we have an
$F\otimes\QQ_\ell$\/-action on $\FGal_\lam$, and we
set $\bV_{\sigma,\lam} = R^1f_*\QQ_\ell
\otimes_{F\otimes\QQ_\ell} \FGal_\lam$, with the
tensor product taken relative to this action and.
Then for $\xi$ as before we set $\bV_{\xi,\lam} =
\otimes_i (\Symm^{k_i}\bV_{\sigma_i,\lam}\otimes
\det\bV_{\sigma_i,\lam})$ with $\sig_i\in\Sigma$
the embedding corresponding to the index $1\leq
i\leq d$. The moduli description gives that
$\bV_{\xi,\lam}$ is pure (in the sense of
\cite{Del2}) of weight $\sum_{i=1}^d (k_i-2
+2\alpha_i) = d\Const$.

\soussection{The dual group}
 \label{sousLgp}
 Let $p$ be a rational prime not dividing $\Disc F$ and
fix some embedding of $\RR$ into $\Qpb$. Then
$\bG(\QQ_p) \simeq \prod_{v|p} \bG_v(\QQ_p)$, where
$\bG_v= \res_{F_v/\QQ_p} \GL_2$. The set $\Sigma$
is a disjoint union $\sqcup_v \Sigma_v$, where
$\Sigma_v$ is the set of the $d_v = [F_v:\QQ_p]$
embeddings of $F_v$ in $\Qpb$. The connected
component $^L\bG_v(\QQ_p)^0$ of the Langlands dual
to $\bG_v(\QQ_p)$ is then the quotient of
$\GL_2^{\Sigma_v}$ by the $\Sigma_v$\/-tuples
$(\dots,z_\sigma,\dots)_{\sigma\in\Sigma_v}$ of
scalar matrices whose product is the identity.
Similarly $^L\bG(\QQ_p)^0$ is the quotient of
$\GL_2^\Sigma$ by the $\Sigma$\/-tuples
$(\dots,z_\sigma,\dots)_{\sigma\in\Sigma}$ of
scalar matrices whose product is the identity. In
particular, $^L\bG(\QQ_p)^0$ is a quotient of
$\prod_{v|p} \,^L\bG_v(\QQ_p)^0$. The Galois group
$\Gal_{\QQ_p} = \Gal(\Qpb/\QQ_p)$ acts naturally on
each $\Sigma_v$ and hence on their union $\Sigma$,
and it is clear that we get an action of
$\Gal_{\QQ_p}$ on $^L\bG_v(\QQ_p)$ and on
$^L\bG(\QQ_p)$. Then $^L\bG_v(\QQ_p)$ is the
semidirect product $^L\bG_v(\QQ_p)^0 \rtimes
\Gal_{\QQ_p}$ and likewise $^L\bG(\QQ_p) =
\,^L\bG(\QQ_p)^0 \rtimes \Gal_{\QQ_p}$. It is clear
that $^L\bG_v(\QQ_p)^0$ acts naturally on
$\CC^{2^{d_v}}$ and that this extends to an
embedding $r_v$ of $^L\bG_v(\QQ_p)$ into
$\GL_{2^{d_v}}(\CC)$. Likewise we get an injective
representation
\[r: \,^L\bG(\QQ_p) \ra \GL_{2^d}(\CC).\]

Now let $\pi_p$ be an irreducible admissible
representation of $\bG(\QQ_p)$ unramified at $p$.
Then $\pi_p$ is a product $\prod_{v|p} \pi_v$ of
unramified principal series representations $\pi_v
\simeq \pi(\mu_{1,v},\mu_{2,v})$ of $\bG_v(\QQ_p) =
\GL_2(F_v)$, with $\mu_{i,v}$ unramified characters
of $F_v^\times$ (see \cite[Ch.\ I.3]{JaLa}). We
will always assume that each
$\pi(\mu_{1,v},\mu_{2,v})$ is infinite-dimensional.
Let $\varpi_v$ be a uniformizer at $v$ and put $a_v
= \mu_{1,v}(\varpi_v)$ and $b_v =
\mu_{2,v}(\varpi_v)$.  Since the central character
of $\pi_v$ maps each $t\in F_v^\times$ to
$t^\Const$, it follows that $|a_vb_v| =
p^{d_v\Const}$. We shall be interested in the size
$|\lam_v|$ of the Hecke eigenvalue $\lam_v = a_v +
b_v$ on $\pi_v$. Recall that $\pi_p$ is {\em
tempered} if and only if all the $\pi_v$, $v|p$
are. This means that $|a_v| = |b_v|$ for all $v|p$.
If this holds --- which is the assertion of the
Ramanujan-Petersson conjecture if $\pi_p$ is a
local component of a cuspidal automorphic
representation --- then $|\lam_v| \leq
2p^{d_v\Const/2}$.

Set $^0t_v = \matr{a_v}{0}{0}{b_v}$. Then the
Satake isomorphism associates to $\pi_v$ the
conjugacy class of the image in $^L\bG_v(\QQ_p)$ of
the element
\[t(\pi_v) = (\,^0t'_v,\Frob_p)\in
\GL_2(\CC)^{d_v} \rtimes \Gal_{\QQ_p},\] where
$^0t'_v = (\,^0t_v, \Id_{2\times 2}, \dots,
\Id_{2\times 2})$. (see e.g. \cite[Exp.\
VI.4]{VSL}). It follows that the Satake isomorphism
associates to $\pi_p$ the element
\[t(\pi_p) = (\,^0t_p,\Frob_p) \in \,^L\bG(\QQ_p)^0 \rtimes
\Gal_{\QQ_p},\] where $^0t_p$ is the image of
$(\dots,\,^0t'_v,\dots)$ in $^L\bG(\QQ_p)^0$.

We now recall the formula for the trace of
$r(t(\pi_p))^m$ for any integer $m\geq 1$. Put $l_v
= \gcd(d_v,m)$, and write $m = q_vd_v + r_v$ with
integers $q_v \geq 0$ and $0\leq r_v < d_v-1$ . We
then have
\[t'(\pi_v)^m =
((\underbrace{^0t_v^{q_v+1},\dots,\,
^0t_v^{q_v+1}}_{r_v\text{\ times}},
\underbrace{^0t_v^{q_v},\dots,\,
^0t_v^{q_v}}_{d_v-r_v \text{\ times}}),
\Frob_p^m).\]
 A straightforward calculation (\cite{Lan} or loc.\
cit.) then gives
\begin{equation}
\label{Trrtp}
 \Tr r(t(\pi_p))^m = \prod_v \Tr r_v(t(\pi_v))^m
  = \prod_v (a_v^{m/l_v} + b^{m/l_v})^{l_v}.
\end{equation}

\soussection{The naive Lefschetz number}
 \label{souslef}
 Write $U = U^p U_p$ with $p$ as before, so that in particular
$U_p\subset \bG(\QQ_p)$ and $U^p \subset
\bG(\AA_\rf^p)$. Let $m\geq 0$ be an integer and
take $a\in \bG(\AA_\rf^p)$. Following \cite{Lan,
BL} one defines functions $f_\xi$ on $\bG(\RR)$,
$\phi_a = \phi_{U^paU^p}$ on $\bG(\AA_\rf^p)$,
$h^m_p$ on $\bG(\QQ_p)$ and $f^G(m,a) = f_\xi
\phi_a^p h_p^m$ on $\bG(\AA) = \bG(\RR)
\bG(\AA_\rf^p) \bG(\QQ_p)$ having the following
properties:
\begin{enumerate}
\item $\Tr\pi_\infty(f_\xi)$ is the multiplicity
of $\xi$ in $\pi_\infty$ for each representation
$\pi_\infty$ of $\bG(\RR)$
\item $\phi_a$ is the characteristic function of $U^paU^p$;
\item $\Tr \pi_p(h^m_p) = p^{md/2} \Tr r(t(\pi_p)^m)$ for
each infinite-dimensional unramified representation
$\pi_p$.
 \end{enumerate}
Let $T(a)$ be the Hecke operator $T(a)$ associated
to $U^paU^p$. Then $T(a) \times \Frob_p^m$ acts as
a correspondence on $\HB_U \times \Fpb$. If $m$ is
sufficiently large then the set of fixed points of
this correspondence is finite, and its graph on
$(\HB_U \times \Fpb)^2$ is transversal to the
diagonal. This enables to define its Lefschetz
number by
\[ \Lef(T(a) \times \Frob_p^m , \HB_U \times \Fpb,
\bV_{\xi,\lam}) = \sum_t \Tr(T(a) \times \Frob_p^m,
(\bV_{\xi,\lam})_t), \]
 the sum taken over the (finite) fixed point set
of $T(a)\times \Frob_p^m$ on $S(\HB_U \times
\Fpb$).

Let $\omega_\xi$ be the central character
$\omega_\xi(t_1,\dots,t_d) = (t_1\dots t_d)^\Const$
of $\xi$, and let $\Ltwo^{{\rm dis}}
(\bG(\QQ)(\bZ(\AA_f)\cap U) \bs\bG(\AA),
\omega_\xi)$ be the space of functions $f$ on
$\bG(\AA)$ satisfying
 \begin{enumerate}
 \item $f(g_\QQ z_U z_\infty g)= \omega_\xi(z_\infty)^{-1}f(g)$
for all $g_\QQ \in \bG(\QQ)$, $z_U \in
\bZ(\AA_f)\cap U$, $z_\infty \in \bZ(\RR)$ and $g
\in \bG(\AA)$.
 \item $f|\det|^{\Const/2}$ is square integrable on
$\bZ(\RR)(\bZ(\AA_f)\cap U)\bG(\QQ)\bs\bG(\AA)$.
\end{enumerate}

As the notation indicates, $\Ltwo^{{\rm dis}}
(\bG(\QQ)(\bZ(\AA_f)\cap U) \bs\bG(\AA),
\omega_\xi)$ is a discrete sum of representations
of $\bG(\AA)$.

Using the moduli interpretation (see \cite{Mil} and
\cite[Expos\'e V]{VSL}), the results of Langlands
(loc.cit.) together with \cite[Section 3.3]{BL}
give the following key

\begin{theorem}
\label{lefpi}
 We have
\[ \Lef(T(a)\times\Frob_p^m, \HB_U\times\Fpb, \bV_{\xi,\lam})
= \sum_\Pi \Tr\Pi(f^G(m,a)),\]
 where the sum is over all irreducible representations $\Pi$ of
$\bG(\AA)$ which occur in the discrete spectrum
$\Ltwo^{{\rm dis}} (\bG(\QQ)(\bZ(\AA_f)\cap U)
\bs\bG(\AA), \omega_\xi)$.
\end{theorem}
(In this formula we are implicity using the strong
multiplicity one theorem for $\GL_2$.)

\soussection{The true Lefschetz number} The problem
now is to identify Langlands's ``naive'' Lefschetz
number as the sum of the local terms in an actual
Lefschetz trace formula. An important technical
issue is the contribution of the boundary. A
conjecture of Deligne asserts it is $0$ for
$T(a)\times\Frob_p^m$, with a given $a$, if $p$ is
sufficiently large. In his thesis Rapoport
constructed smooth (toroidal) compactifications of
$\HB_U$ over $\ZZ[1/ND]$ in which the complement of
$\HB_U$ is a relative normal crossing divisor
\cite[Corollaire 5.3]{Rap}. For all but finitely
many primes $p$ there also exists the Baily-Borel
compactification, whose boundary is the finite set
of cusps, and the toroidal compactification is a
blow-up of it. At such primes $p$, Brylinski and
Labesse could prove Deligne's conjecture
(\cite[Theorem~2.3.3]{BL}). According to C.-L. Chai
(private communication) this holds for all $p$
prime to $ND$. Alternatively we can use the results
of \cite{Pin}, since Condition~7.2.1 there, which
suffices for Deligne's conjecture, holds for these
compactifications, provided we know also that the
monodromy of our sheaves around the toroidal
resolutions of the cusps are tame. This is manifest
from the higher dimensional Mumford-Raynaud-Tate
parametrization near the cusps \cite[Th\'eor\`me
5.1 and 4.11]{Rap}. The idea is that the
$\ell$\/-adic Tate module of a universal abelian
variety in a punctured $p$\/-adic analytic
neighborhood of a cusp admits a canonical exact
sequence having an \'etale quotient and a
multiplicative sub-object (which are moreover
Cartier dual to one another). The $\ell$\/-adic
monodromy is therefore manifestly tame at $p$.
Finally, we could appeal to \cite{Fuj}, where
Deligne's conjecture is proved in general. Either
way we get the following
\begin{corollary}
\label{leffrob}
 For every $a$ as above there exists an integer $m_0(a)$
such that for each integer $m\geq m_0(a)$ we have
\begin{multline*}
\Lef(T(a) \times \Frob_p^m , \HB_U \times \Fpb,
\bV_{\xi,\lam}) \\
 = \Tr(T(a)\times \Frob_p^m |
 H^*_c(\HB_U \times \Fpb, \bV_{\xi,\lam}),
 \end{multline*}
where as usual $\Tr(\cdot| H^*_c)$ means the
alternating sum $\sum_{i=0}^{2d} (-1)^i \Tr(\cdot |
H^i_c)$.
\end{corollary}

\soussection{Cuspidality and compact support}
Assume that $\pi = \otimes_p \pi_p$ is an
irreducible {\em cuspidal} representation of
$\bG(\AA)$. In particular each $\pi_p$ (and
$\pi_v$) is infinite dimensional, and $\pi\otimes
|\det|^{-\Const/2}$ is unitary. Then $\pi$
contributes to $H^*(\HB_U,\bV_\xi)$ if and only if
$\pi$ has a $U$\/-invariant vector and if
$\pi_\infty$ is isomorphic to
\[ \prod_{i=1}^d
    \pi(|{\bf \cdot}|^{(1-k_i)/2},
          |{\bf \cdot}|^{(k_i-1)/2}\sign^{k_i})\
            \otimes {\det\,}^\Const.\]
We will then say that $\pi$ is of {\em type
$\xi$}\/. In this case $\pi$ contributes to
$H^i(\HB_U,\bV_\xi)$ if and only if $i=d$. In fact
taking the inverse limit over $U$ we set
\[ U^d(\pi_f,\xi) =
     \Hom_{\bG(\AA_f)} (\pi_f, \invlim_U H^d(\HB_U,\bV_\xi)).\]
Then $U^d(\pi_f,\xi)$ is a $2^d$\/-dimensional
space, (see \cite[Section 3.4]{BL}, and notice that
by the strong multiplicity $1$ theorem for $\GL_2$
we do not have multiplicities in our case). The
$\ell$\/-adic analog likewise gives a
$2^d$\/-dimensional $\Gal_\QQ = \Gal(\ov{\QQ}/\QQ)$
representation $U^d_\lam(\pi_f,\xi)$. We now need
the following
\begin{proposition}
If $\pi$ is a cusp form as above, then
$U^d(\pi_f,\xi)$ comes from the compactly supported
cohomology; in other words, the natural map of
forgetting supports
\[
\Hom_{\bG(\AA_f)} (\pi_f, \invlim_U
H^d_c(\HB_U,\bV_\xi))
       \ra U^d(\pi_f,\xi)\]
is an isomorphism.
\end{proposition}
\begin{pf}
See \cite[Cor. 5.5]{Bor} (and also the last comment
in \cite{Har1}). Here is the general idea of the
proof. For any cusp $c$ of $\HB_U$ let $\HB_{U,c}$
be the boundary component corresponding to $c$ in
the Borel-Serre compactification of $\HB_U$. The
cuspidality of $\pi$ implies that the (de-Rham)
cohomology classes of the vector-valued
differential forms $\set{\omega}$ associated to it
restrict to $0$ on each $\HB_{U,c}$. This can be
seen also by verifying that the periods
$\int_\beta\omega$ of these differential forms
around each (vector-valued Borel-Moore)
$d$\/-dimensional cycle $\beta$ is $0$. In fact, we
can move $\beta$ towards the cusp; then on the one
hand the period stays constant, but on the other
hand the cuspidality of $\pi$ makes the $\omega$
and with it the period decrease exponentially. This
forces the period to be $0$, and we conclude that
our $\set{\omega}$\/'s are cohomologous to
compactly supported forms.
\end{pf}

As a corollary we get that $U^d(\pi_f,\xi)$, which
initially was derived from $H^d(\HB_U,\bV_\xi)$, is
in the image of the ``forget support'' map from
$H^d_c(\HB_U,\bV_\xi)$. This image is the {\em
parabolic} cohomology
$\tilde{H}\vphantom{H}^d(\HB_U,\bV_\xi)$, and it
exists in $\ell$\/-adic cohomology. As the sheaf
$\bV_{\xi,\lam}$ is pure of weight $d\Const$, a
weight argument (see \cite{Del2}) gives that
$\tilde{H}\vphantom{H}^d(\HB_U,\bV_\xi)$ is pure of
weight $d(\Const+1)$.  The same is therefore true
for $U^d_\lam(\pi_f,\xi)$.

\soussection{The Ramanujan-Petersson Conjecture} We
partially recall our previous notation. With $\xi$
as before (in particular, of central character
$\omega_\xi = \Nm^\Const$), let $\pi=\prod_p\pi_p$
be an irreducible cuspidal representation of
$\bG(\AA)$ of infinity type $\xi$, so that
$\pi_\infty$ is the discrete series representation
corresponding to $\xi$ as above. Let $p$ be a
rational prime which is prime to the level of $\pi$
and to $\Disc F$. For a prime $v$ of $F$ above $p$
let $d_v$ denote the degree of $F_v/\QQ_p$, and let
$T_v$ be the Hecke operator at $v$, defined as the
action of the double class
$U_v\matr{\varpi_v}{0}{0}{1}U_v$, with $U_v \simeq
\GL_2(\cO_{F,v})$ having measure $1$. Here
$\varpi_v$ is a uniformizer at $v$. The Satake
parameters of $\pi_v$ are given as before in terms
of $a_v$ and $b_v$ (up to order). The eigenvalue of
$T_v$ on $\pi$ (or on $\pi_v$) is $\lam_v = a_v +
b_v$. We now obtain the following version of the
Ramanujan-Petersson Conjecture:
\begin{theorem}
\label{RP} The representations $\pi_p$ and $\pi_v$
for $v|p$ are tempered. More precisely, each of
$a_v$, $b_v$ is a Weil number of weight
$d_v+\Const$, i.e. $a_v$ and $b_v$ are algebraic
numbers such that under any embedding into $\CC$ we
have
 \[ |a_v| = |b_v| = p^{d_v \Const/2}.\]
Consequently, $|\lam_v| \leq 2p^{d_v \Const/2}$.
\end{theorem}

\begin{pf} Using our previous notation, there exists a
finite linear combination $T = \sum c_i T(g_i)$ of
Hecke operators away from $p$ which acts as
projection onto the $U^p$ invariants of the
$\pi_f^p$\/-isotypical part. By the strong
multiplicity one theorem this is a projection onto
the $U^p$ invariants of the $\pi_f$ isotypical
part. Let $m_0$ be an integer such that Deligne's
conjecture holds for $T(g_i) \times \Frob_p^m$ on
$\HB_U \times \Fpb$ for all $i$ and for all $m \geq
m_0$. From Theorem~\ref{lefpi} and
Corollary~\ref{leffrob} applied to
$T\times\Frob_p^m$, with $T$ viewed as the above
linear combination, and from the definition of
$U^d_\lam(\pi_f,\xi)$ we get
\[ \Tr(\Frob_p^m|U^d_\lam(\pi_f,\xi)) = \Tr\pi_p(h_p^m).\]
The purity of the parabolic cohomology and the
property of $h_p^m$ from
Section~\ref{secHMF}.\ref{souslef} give the
existence of $2^d = \dim U^d_\lam(\pi_f,\xi)$ Weil
numbers $\set{w_i}_i$ of weight $d(\Const+1)$ so
that for every large enough $m$ we have
\[ \sum_i w_i^m = p^{md}\Tr r(t(\pi_p)^m) = p^{md}\prod_{v|p}
(a_v^{m/l_v} + b_v^{m/l_v})^{l_v},\]
 by equation~(\ref{Trrtp}), where as before $l_v =
\gcd(m,d_v)$. As in \cite[Theorem~3.4.6]{BL} this
implies that each $\pi_v$ is tempered, namely that
$|a_v| = |b_v|$ for every embedding of $a_v$, $b_v$
into $\CC$. That the weight is as it should be then
follows from our knowledge of the product
$|a_vb_v|$, as was explained in
Section~\ref{secHMF}.\ref{sousLgp}.
\end{pf}

Theorem~\ref{RPp} follows when one makes the
necessary normalizations. In particular one takes
$\Const = \max_i k_i$.

\section{Special cases}
\label{secspec}
 We now give two special cases when the
conditions \ref{specase} are satisfied. The first
one was accessible in \cite{JL8}, while the second
one was not.

\noindent {\bf A.}\ Take $F = \QQ$ and let $B$ be a
quaternion algebra over $\QQ$ of discriminant $2$
(it is unique up to an isomorphism). A model for
$B$ is given by the usual rational quaternion
algebra generated over $\QQ$ by $\qi$, $\qj$
satisfying $\qi^2 = \qj^2 = -1$ and $\qi\qj = -
\qj\qi$. It is very easy to see that
Conditions~\ref{specase} are satisfied with $N_0
=2\ZZ$. Taking $g = 1$ prime $p$ satisfying
$p\equiv 1\mod{4}$, one gets the
Lubotzky-Phillips-Sarnak graphs \cite{LPS}. These
are Cayley graphs on $\PSL(2,\FF_N)$, for any prime
$N\equiv 1\mod 4$ which is different from $p$ and a
a square modulo $p$. The set of generators is the
reductions modulo $N$ of the $p+1$ norm $p$
integral quaternions $\equiv 1\mod{2}$.

The case of $g = 2$ distinct primes $p$, $q$
satisfying $p\equiv q \equiv 1 \mod{4}$ appeared in
a different context in \cite{Moz}. The resulting
two (blue and red) communication networks are
realized as the same Cayley graphs as before on
$\PSL(2,\FF_N)$, for any prime $N\equiv 1\mod 4$
which is different from $p$ and $q$ and is a a
square modulo both. The sets of generators are the
reductions modulo $N$ of the $p+1$ (respectively
$q+1$) norm $p$ (respectively norm $q$) integral
quaternions $\equiv 1\mod{2}$. They have the square
property. Thus each node has $p+1$ blue neighbors
and $q+1$ red ones.

\noindent {\bf B.}\ We shall now construct a blue
and a red Ramanujan communication networks having
the square property over the same set of vertices
$\PSL(2,\FF_N)$ in which each node has $p+1$ blue
and $p+1$ red neighbors. Here $p$ is any prime
$\equiv 1, 9\mod{20}$ and $N\neq p$ is a prime
satisfying $N\equiv 1,9\mod{20}$ and which is a
square modulo $p$. There are other conditions which
$N$ must satisfy. These cannot be stated as simply,
but we will see that they are satisfied by a set of
primes $N$ of positive density which can be
explicitly described.

Take $F = \QQ(\sqrt{5}) \subset \RR$. Let $\infty_1
= \Id$, $\infty_2$ be the real primes of $F$. The
fundamental unit of $F$ is $\tau = (1 +
\sqrt{5}\,)/2$. We will need the following lemma
from Class Field Theory:
\begin{lemma}
\label{CFT}
 1.\ Let $N$ be a rational prime such that $\leg{-1}{N}
= \leg{5}{N} = 1$ (in other words $N\equiv 1,9
\mod{20}$). It is then possible to write $N = a^2 -
20b^2$ with $a,b$ integers and $a>0$ odd. Then
$\tau$, considered modulo $N$ via either choice of
$\sqrt{5} \mod{N}$, is a square modulo $N$ if and
only if $a + 2b \equiv 1\mod{4}$, and this holds
for a set of primes $N$ of Dirichlet density $1/8$.

\noindent 2.\ Let $p\neq N$ be a rational prime
satisfying with $\leg{5}{p} = \leg{-1}{p} = 1$ and
use part~1.\ to write $p = a^2 - 20b^2$ for
integers $a$, $b$. Choose a $\sqrt{5}$ modulo $N$.
Then the set of primes $N$ satisfying the criterion
in part~1.\ and for which both $a\pm 2b\sqrt{5}$
are squares modulo $N$ has Dirichlet density
$1/32$. Such primes $N$ satisfy in particular
$\leg{p}{N} = 1$.
\end{lemma}
\begin{pf}
If $5$ is a square modulo $N$ then $N$ splits in
$F$, and since the class number of $F$ is $1$ we
get $\pm N = \nu\nub$, with $\nu$ and $\nub$
conjugate integers of $F$. Since $-1$ is the norm
of $\tau$ and since $\tau$ reduces modulo $2$
(which is inert in $F$) to a generator of
$\FF_4^\times$, with $\FF_4 = \cO_F/2\cO_F$, we can
assume that the sign is $+$ and that $\nu \equiv 1
\mod{2\cO_F}$. Writing $\nu = a + b'\sqrt{5}$, we
get that $a$ and $b'$ are integers with $a$ odd and
$b' = 2b$ even. Multiplying by $-1$ we may also
assume that $a> 0$. A choice of a square root of
$5$ modulo $N$ then determines the sign of $b$, or
equivalently the prime $\nu = a + 2b\sqrt{5}$ above
$N$. We then ask when is the ideal $\nu\cO_F$ split
in $L= F(\sqrt{\tau})$. Since $N\equiv 1\mod{4}$,
it splits in $\QQ(\sqrt{-1}\,)$, so that $\nu$ (and
$\ov{\nu}$\/) split in $F(\sqrt{-1}\,)$. Therefore
$\nu$ splits in $F(\sqrt{\tau}\,)$ if and only if
$N$ is completely split in
$\QQ(\sqrt{5},\sqrt{-1},\sqrt{\tau}\,)$, which is
Galois of degree $8$ over $\QQ$. By the
\v{C}ebotarev density theorem, this happens for a
set of primes $N$ of Dirichlet density $1/8$
(namely, for half the primes $\equiv 1,9
\mod{20}$). It remains to determine the condition
for $\nu$ to split in $F(\sqrt{\tau}\,)$.

The Id\`ele class group character $\chi$ of $F$
which cuts $L$ is of order $2$ and unramified
outside $2$. Since $\infty_1(\tau) > 0$ and
$\infty_2(\tau) < 0$, we see that $\chi$ is trivial
on $F^\times_{\infty_1}$ but not on
$F^\times_{\infty_2}$. For a prime $v$ of $F$ set
$U_v = \cO_{F,v}^\times$ and let $U_{v,n}$ be its
subgroup of elements congruent to $1$ modulo $v^n$
for any $n\geq 1$. Then $\chi$ must be trivial on
the product $U_\infty$ of the connected components
of $F_{\infty_i}^\times$ and on $U = \prod_{\neq 2}
U_v$. The units $\cO_F^\times = \pm\tau^\ZZ$
surject onto $U_2/U_{2,1} \times \prod_i
\pi_0(F_{\infty_i})$, and the kernel of this
surjection is $\tau^{6\ZZ}$. Since the class number
of $F$ is $1$, we get that $\chi$ is determined by
its values on
\[ C =
F^\times \bs \AA_F^\times /U U_\infty \simeq
\tau^{6\ZZ} \bs U_{2,1} /(U_{2,1})^2 \,.\]

First we calculate $C$. Through the $2$\/-adic
logarithm map, $U_{2,1} \simeq \mu(F) \oplus L$,
where $L= 2\ZZ_2 + 4\cO_{F,2}$ and $\mu(F) =
\set{\pm 1}$ is the group of the roots of unity of
$F$.
Since $\tau^6 = 9 + 4\sqrt{5} \equiv 5\mod{8\cO_F}$
we get that
\[C \simeq \mu(F) \oplus 2\cO_F/(2\ZZ + 4\cO_F)
\simeq (\ZZ/2\ZZ)^2 \,,\] with generators $\ov{-1}$
and $2\ov{\tau}$ for the respective two factors.

To determine $\chi$ we use the product formula
$\prod_v \chi_v(x) = 1$ for any $x\in F^\times$,
the product taken over all places of $F$. For
$x=-1$ we get $\chi(\ov{-1}) = \chi_2(-1) =
                      \chi_{\infty_2}(-1) = -1$,
and since $\tau^3 = 2 + \sqrt{5} = 1+2\tau$, we
also get
\[ \chi(2\ov{\tau}) = \chi_2(1 + 2\tau)
    = \chi_2(\tau)^3  = -1\,.\]

In particular, for $x,y\in \ZZ_2$ we have $\chi_2(1
+ 2(x + y\tau)) = (-1)^{x + y}$. It follows the
prime $\nu\cO_F$ splits in $F(\sqrt{\tau}\,)$ if
and only if $\chi_\nu(\nu) = -1$. But
\[\chi_\nu(\nu) =
\chi_2(\nu)\chi_{\infty_2}(\nu)
 = (-1)^{(a-1)/2 + b} \cdot 1
   = (-1)^{(a-1)/2 + b}\,.\]
As claimed. (We write $\chi_\nu$ etc.\ for the
local component of $\chi$ at the place defined by
the prime $\nu\cO_F$\/.)

\noindent 2.\ Write $p = \pi \pib$, with (say) $\pi
= a + 2\sqrt{b}$. Observe that the field
$F(\sqrt{\pi},\sqrt{\pib})$ is Galois over $\QQ$
and contains $\sqrt{p}$. As an extension of $F$ it
is ramified only over $\pi$ and $\pib$: this is
because $(1 + \sqrt{\pi})/2$ is an algebraic
integer, as its trace and norm to $F$ are $1$ and
$(1-a)/2 - b\sqrt{5}$. Hence it is linearly
disjoint from $F(\sqrt{-1},\sqrt{\tau})$ over $F$.
Applying the \v{C}ebotarev density theorem to the
degree $32$ extension
\[\QQ(\sqrt{-1},\sqrt{5},\sqrt{p},\sqrt{\tau},
     \sqrt{\pi})\]
of $\QQ$ gives the required density statement,
concluding the proof of the lemma.
\end{pf}

Now set $B_F = B \otimes F$, with $B$ the rational
quaternions as before. For the facts we need about
$B_F$ see \cite[Chapter 5]{Vig}. In particular,
$B_F$ is ramified precisely at the two infinite
primes of $F$.

The class number of $B_F$ is $1$, and all maximal
orders in $B_F$ are conjugate. To describe the sets
of generators of our red and blue networks we need
to explicitly compute in one maximal order $\cM$.
We take for $\cM$ the $\cO_F$\/-submodule of $B_F$
generated by
\begin{equation}
\begin{array}{rcl}
e_1 & = & (1 + \tau^{-1}\qi + \tau\qj)/2 \\
e_2 & = & (\tau^{-1}\qi + \qj + \tau\qk)/2 \\
e_3 & = & (\tau\qi + \tau^{-1}\qj + \qk)/2 \\
e_4 & = & (\qi + \tau\qj + \tau^{-1}\qk)/2\
\end{array}
\end{equation}
We have  $\cM/2\cM \simeq \Matwo(\FF_4)$, with
$\FF_4 = \cO_F/2\cO_F$ a field with $4$ elements.
The reduction map induces an exact sequence
\[1 \ra (\cO_F + 2\cM)^\times \ra \cM^\times \ra
\PGL(2,\FF_4) \ra 1 \,.\]
 From our choice of $\cM$ it follows that
$\cO_F + 2 \cM$ is the $(\cO_F$\/-) sub-order of
$\cM$ given by
\[ \set{a + b\qi + c\qj +d\qk \mid | a,b,c,d\in \cO_F,
   \ b + \tau c + \tau^{-1}d \in 2\cO_F}\,.\]
Since $\tau$ reduces modulo $2$ to a generator (of
order $3$) of $\FF_4 ^\times$, we obtain the
following
\begin{lemma}
Every element in $\cM$ whose norm to $\cO_F$ is
odd, namely not in $2\cO_F$, can be uniquely
written as $u(a + b\qi + c\qj + d\qk)$, with $u \in
\cM^\times$ and $a,b,c,d\in \cO_F$ satisfying
\begin{enumerate}
\item  $b + \tau c + \tau^{-1}d \in 2\cO_F$;
\item $a \equiv 1 \mod{2\cO_F}$;
\item $1\leq a < \tau^3$.
\end{enumerate}
When $u = 1$ we shall say that our element is in
{\em normal form}.
\end{lemma}

As a corollary we see that Conditions~\ref{specase}
hold with $N_0 = 2\cO_F$.

Now write $4p = a^2 - 20b^2$, with $a$, $b$
integers (say, minimal positive), necessarily of
the same parity. The elements $\pi = (a +
2b\sqrt{5}\,)/2$ and $\pib = (a - 2b\sqrt{5}\,)/2$
are the primes above $p$ in $\cO_F$. Let $\cMF
=\cO_F[\qi,\qj]$ be the order of $\cO_F$\/-integral
quaternions. We shall say that an element $x\in
\cMF$ is in normal form if $x \equiv 1 \mod{2}$ and
$\Tr x = \Tr_{B_F/F} x$ is totally positive (namely
positive for both real embeddings) and satisfies
$\tau^{-3} < \Tr x \leq \tau^3$. We then have the
following
\begin{lemma} {\rm 1}.\ There are precisely $p+1$
elements in normal form $\gamma_i$ (respectively
$\gammab_i$), $1\leq i \leq p+1$, in $\cMF$ whose
norm is $\pi$ (respectively $\pib$).

\noindent {\rm 2}.\ For any indices $1\leq i, j\leq
p+1$ there exist unique indices $1\leq i', j' \leq
p+1$ and a unique unit $u \in \cMF^\times =
              \langle \pm \pi^{3\ZZ} p^\ZZ\rangle $
such that
\[\gamma_i \gammab_j = \gammab_{j'} \gamma_{i'} u\,.\]
\end{lemma}
\begin{pf}
This is merely an explicit restatement of
Proposition~\ref{generators}.
\end{pf}

Now let $N$ be a rational prime congruent to $1$ or
$9$ modulo $20$. Choosing square roots of $-1$ and
of $5$ in the prime field $\FF_N$ allows us to map
$\cM$ homomorphicallly onto $\Matwo(\FF_N)$ by
\[ \sqrt{5} \mapsto
   \matr{\sqrt{5}}{0}{0}{\sqrt{5}}\,,\;\;
   \qi \mapsto \matr{\sqrt{-1}}{0}{0}{-\sqrt{-1}}
   \,,\;\;\text{and}\;\;
   \qj \mapsto \matr{0}{-1}{1}{0}\,.\]
Proposition~\ref{explicit} now gives us that we
have constructed two Ramanujan communication
networks as the two Cayley graphs on
$\PSL(2,\FF_N)$ having the reductions modulo $N$ of
the $\gamma_i$\/'s (respectively the
$\gammab_j$\/'s) for generators.

\appendix

And this is my Appendix.




\end{document}